\documentclass[12pt]{article}
\makeatletter
\def\section{\@startsection{section}{1}{\z@}{-3.5ex plus -1ex minus -2.ex}
{2.3ex plus .2ex}{\Large\bf}}

\def\subsection{\@startsection{subsection}{2}{\z@}{-3.25ex plus
 -1ex minus -2.ex}
{1.5ex plus .2ex}{\bf}}

\def\vsn{\vskip 1pc \noindent}
\def\f{\newline}
\def\e{\varepsilon}

\def\cost{{\rm cost}}

\def\rr{{\bf R}}
\def\nn{{\bf N}}

\def\k{\kappa}

\newcommand{\be} {\begin{equation}}
\newcommand{\ee} {\end{equation}}
\newcommand{\bd} {\begin{displaymath}}
\newcommand{\ed} {\end{displaymath}}

\newcommand{\bq}{\begin{eqnarray}}
\newcommand{\eq}{\end{eqnarray}}
\newcommand{\bqn}{\begin{eqnarray*}}
\newcommand{\eqn}{\end{eqnarray*}}
\newcommand{\ba}[1]{\begin{array}{#1}}
\newcommand{\eqa}{\end{array}}
\def\qed{
   \\[-4ex]
  \hbox to \hsize{\hfill \vrule height 1.6ex width 1.5ex
  depth -.1ex}}

\begin{document}

\bibliographystyle{alpha}

\begin{center} {\Large {\bf
Adaptive mesh selection asymptotically 
guarantees a prescribed local error for systems of initial value problems
}\footnotemark[1] }
\end{center}
\footnotetext[1]{ ~\noindent This research was partly supported by the Polish
NCN grant - decision No. DEC-2017/25/B/ST1/00945    and by the Polish   Ministry of Science  and Higher Education
\vsn}               
  
\medskip 
\begin{center}
{\large {\bf Boles\l aw Kacewicz \footnotemark[2]  }}
\end{center}
\footnotetext[2]{ 
\begin{minipage}[t]{16cm} 
 \noindent
{\it AGH University of Science 
and Technology, Faculty of Applied Mathematics,\\
\noindent  Al. Mickiewicza 30, paw. A3/A4, III p., 
pok. 301,\\
 30-059 Cracow, Poland 
\newline
E-mail B. Kacewicz: $\;\;$ kacewicz@agh.edu.pl }  
\end{minipage} }

\thispagestyle{empty}
$~$
\vsn
\vsn
\begin{center} {\bf{\Large Abstract}} \end{center}
\noindent
We study  potential advantages of adaptive mesh point selection for the solution of systems of initial value problems. 
 For an optimal order discretization method, we
propose an algorithm for successive selection of the mesh points, which only requires evaluations of the 
right-hand side function.  The selection (asymptotically) guarantees  that the  maximum local error of the method does not exceed a prescribed level. 
The usage of the algorithm is not restricted to the chosen method; it can also be applied  with any method from
 a general class.
 We  provide a rigorous analysis of the cost of the proposed algorithm. It is shown that the cost is almost minimal,   up to absolute constants,
among all mesh selection  algorithms.   For illustration, we  specify the advantage of the adaptive mesh over 
 the uniform one.  Efficiency of the adaptive algorithm   results from automatic adjustment of the successive mesh points
 to the local behavior of the solution.
  Some numerical results  illustrating  theoretical findings are reported. 
\newpage
\noindent
\section {\bf{\Large Introduction}} 
\noindent
We deal with the solution of systems of initial value problems (IVPs)
\be
z'(t)=f(t, z(t)),\;\; t\in [a,b],\;\;\;\; z(a)=\eta,
\label{1}
\ee
where $a<b$,   $\eta\in \rr^d$ and  $f: [a,b]\times \rr^d \to \rr^d$ is a $C^r$ function.
We  study   how much adaptive  mesh points improve efficiency of
algorithms for solving (\ref{1}).  For many years, adaption has been  a standard tool 
in numerical packages.   Well known examples include the package QUADPACK \cite{Pie} for numerical integration, or,
 among other solvers,  the DIFSUB procedure by C.W. Gear or the library ODEPACK by A. Hindmarsh  for IVPs.
Many authors  have  reported   superiority of adaption over nonadaption for (\ref{1}),  based on  numerical results
 for a number of computational examples.   One can cite as an example  papers   such as  \cite{jackiewicz},
\cite{lyness} or \cite{mazzia}.
Practical efficiency gives us   considerable knowledge about the power  of adaption. Conclusions are 
however  not complete;    an analysis of theoretical aspects is missing in  many cases of step size control strategies. 
 For instance, most often,  step size control devices are not supported by   cost analysis. 
The use of variable step size not only improves the efficiency of methods for solving regular problems (\ref{1}), 
but it also allows us   to manage singularities.  Considerable progress has been made in  rigorous
analysis of adaption in that case,  see e.g.  \cite{KP}, \cite{PlaWas}.      Adaption allows us in many cases
to maintain for singular problems the order of convergence known from the regular case.
\f
Advantages of adaption  for regular problems are of different type.  Integration of scalar $C^4$ functions by the Simpson rule
  have  been  recently studied    in \cite{plaskota}.
It is shown  in \cite{plaskota} that 
adaption does not improve the order of convergence, but it can reduce the asymptotic constant of the method.  In the similar spirit,  adaption 
has been considered for univariate approximation and minimization in \cite{hickernell}.  
 For scalar autonomous problems (\ref{1}), adaptive  mesh selection
 has been recently studied in  \cite{NumAlg}. An adaptive strategy has been proposed and the cost analyzed,  based on specific properties of scalar autonomous equations.   The particular technique used  does not allow us to extend the results  from \cite{NumAlg}  to systems of IVPs.
\f
  In the present paper,  we consider general systems (\ref{1}), with 
a $C^r$ right-hand side function $f$.  Our contribution can  be summarized as follows:
\begin{itemize}
\item  For a maximal order method for solving the system (\ref{1}), we propose a new mesh selection algorithm that guarantees  the local error
of the method at a precribed level $\e$, for sufficiently small $\e\in (0,1)$. At each time step, we select the step size
and compute an approximation to the solution, which requires two runs of the basic method.
 Information  about the function $f$ only consists of function evaluations.
We show that the mesh selection algorithm can be applied to a general class of methods for solving (\ref{1}).
\item We rigorously analyze the cost of the method equipped  with the mesh selection procedure. We show that the cost is  minimal 
among all mesh selection strategies, provided that we accept some absolute constants. 
\item We specify in a quantitative  way an advantage of the adaptive mesh over  the uniform one. 
\end{itemize}
The adaptive mesh selection for  systems (\ref{1}) does not improve the speed of convergence of algorithms.
A potential gain of adaption   lies in  reducing a coefficient in the error bound.
\f
The paper is organized as follows. 
 We formulate the problem in Section 2. An algorithm $\phi^*$ with maximal order 
is defined in Section 3, and the error analysis is given in Theorem~ 1. In Section 4, we discuss upper bounds on the local error of $\phi^*$,
in particular, we  show a constructive upper bound in Theorem 2. Section 5 contains  the main algorithm ADAPT-MESH which combines $\phi^*$
with a new  mesh selection algorithm.  A generalization 
of the algorithm ADAPT-MESH is given in Section 6. Section 7 contains  a cost analysis of  the algorithm ADAPT-MESH
compared to other  algorithms, see Theorem 3.  Possible advantage of the 
adaptive mesh over  the uniform one is discussed in detail. 
   Finally, some  results of  numerical experiments are reported in Section 8.
\vsn
\section {\bf{\Large Problem formulation}} 
\noindent
We consider problem (\ref{1}) with a continuous  function $f: [a,b]\times \rr^d \to \rr^d$ 
such that for some $L\geq 0$ 
\be
\| f(t,y) -f(t,\bar y)\| \leq L\|y-\bar y\| \;\;\;\mbox{ for } t\in [a,b],\;\; y,\bar y\in \rr^d. 
\label{lip-ukl}
\ee
Here and in what follows $\|\cdot\|$ denotes 
the maximum norm in $\rr^d$.  It follows that there exists a unique solution $z$ of (\ref{1}) defined on $[a,b]$. 
Let $r\geq 1$. We assume that
$f$ is a regular function in a subset of its domain, $f\in C^r([a,b]\times D)$,  where
\be
D=\{y\in \rr^d:\;  \|y\| \leq \sup\limits_{t\in [a,b]} \|z(t)\| +1\}.
\label{D}
\ee
The class of functions $f$ satisfying the above assumptions will be denoted by $F_r$.
\f
 Let $m\in \nn$. We wish to compute an approximation to the solution $z$  in $[a,b]$.    For 
 $m+1$ mesh points $a=x_{0,m}< x_{1,m}<\ldots < x_{m,m}=b$, we do it by computing approximations
$l_i$ to $z$ in the subintervals $[x_{i,m}, x_{i+1,m}]$, $i=0,1,\ldots,m-1$.
\f 
Let ${\ell}(m)$ be any nonincreasing sequence convergent to $0$ as $m\to \infty$.
We consider for any $f$  a class of partitions of $[a,b]$ defined as follows. We assume that 
    there exist $K=K(f,a,b,\eta)$ and $k_0=k_0(f,a,b,\eta)$ such that for all $m\geq k_0$ and any partition it holds
\be
\max\limits_{0\leq i\leq m-1}(x_{i+1,m} - x_{i,m}) \leq K\, {\ell}(m).
\label{podzial}
\ee
Thus, the partitions under consideration are uniformly normal. 
Note that we always have $\max\limits_{0\leq i\leq m-1}(x_{i+1,m} - x_{i,m})\geq (b-a)/m$ for any $m\geq 1$. Thus, the condition 
 (\ref{podzial}) implies that ${\ell}(m)$ cannot go to zero faster than $1/m$. Note that the convergence of ${\ell}(m)$ to zero can be arbitrarily slow, and 
the constant $K$ can be arbitrarily large.
We shall omit in the sequel the second
subscript $m$, remembering that the choice of  points $x_i$  can be different for varying $m$.
\f
For  a given $y_i\in \rr^d$, we denote by $z_i$  the solution of a local problem
\be
z_i'(t)=f(t, z_i(t)),\;\; t\in [x_i,x_{i+1}],\;\;\;\; z_i(x_i)=y_i.
\label{2}
\ee
If $l_i$ is an approximation to $z_i$ given by a certain method, then  
 local errors of the method are given by 
$$
\sup\limits_{t\in [x_i,x_{i+1}]} \|z_i(t)-l_i(t)\|, \;\;\;\;i=0,1,\ldots,m-1 .
$$
Our aim is to 
select possibly small $m$ and  mesh points $\{x_i\}_{i=0}^m$  such that the local errors remain at a precribed level  $\e>0$.

\section {\bf{\Large The basic method $\phi^*$ and its error}} 
\noindent
The basic method makes use of the approximate Picard iteration, an idea that turned out useful in  several contexts, 
see e.g. \cite{daun},  \cite{KP}.  Let $m\in \nn$ and $x_0=a<x_1<\ldots< x_m=b$ be mesh points satisfying (\ref{podzial}).
Let $y_0=\eta$. For a given $y_i$, we define approximations $l_{i,j}$ in $[x_i,x_{i+1}]$ as follows.
\f
 We set 
$l_{i,0}(t)\equiv y_i$.
Let $l_{i,j}$ be given.  Denote by $t_0,t_1,\ldots, t_{r-1}$ the equidistant nodes in $[x_i, x_{i+1}]$, with 
$t_0=x_i$  for $r=1$ and   $t_0=x_i$, $t_{r-1}= x_{i+1}$ for $r\geq 2$.
(The points $t_k$ depend on $i$;   we shall omit this index  in the notation.) 
We define $q_{i,j}$ to be the Lagrange interpolation polynomial of degree  at most $r-1$ for the function
$g_{i,j}(t)=f(t,l_{i,j}(t))$ based on the nodes   $t_0,t_1,\ldots, t_{r-1}$,
\be
q_{i,j}(t)=\sum\limits_{k=0}^{r-1} g_{i,j}(t_k) \prod\limits_{p=0,p\ne k}^{r-1} \frac{ t-t_p}{ t_k-t_p}, \;\;\;\; t\in [x_i,x_{i+1}],
\label{lagrpol}
\ee
where $\prod\limits_{p=0, p\ne k}^0 =1$. An approximation $l_{i,j+1}$ is given by
\be
l_{i,j+1}(t)=y_i+\int\limits_{x_i}^t q_{i,j}(\xi)\, d\xi,\;\;\;  t\in [x_i,x_{i+1}].
\label{lposr}
\ee
The final approximation in $[x_i,x_{i+1}]$ is given by $l_{i,r+1}$.  To complete the definition, we set $y_{i+1}=l_{i,r+1}(x_{i+1})$.
\f
For  $t\in [a,b]$ we define a continuous  approximation to $z$ by
\be
l_{r+1}(t)= l_{i,r+1}(t),\;\;\; t\in [x_i,x_{i+1}].
\label{lost}
\ee
The transformation that assignes to $f$ the approximation $l_{r+1}$ will be denoted by $\phi^*$, $\phi^*(f)(t)=l_{r+1}(t)$, $t\in [a,b].$
\f
The following theorem provides  error analysis of the method  $\phi^*$.  The proof follows usual lines of the analysis of approximate 
Picard iteration,  it  is  however focused on our specific  requirements.   We shall need in the next sections the error bound (\ref{prop1}) 
for a non-uniform mesh,  as well as specific local error bounds derived in the body of the proof.
   Let $h_i=x_{i+1}-x_i$.
\vsn
{\bf  Theorem 1}$\;\;$ {\it  Let $f\in F_r$. There exists $m_0$ such that for all $m\geq m_0$ and any $\{x_i\}_{i=0}^m$ satisfying
 (\ref{podzial}), the global error of $\phi^*$ at $f$ satisfies
\be
\sup\limits_{t\in [a,b]} \|z(t)-l_{r+1}(t)\| \leq M\max\limits_{0\leq j\leq m-1} h_j^r,
\label{prop1}
\ee
where $M=   \exp(L(b-a))(b-a)\left( 2D_r /r! + 1/2 \right)$, and  the number $D_r$, only dependent on $f$, $r$, $a$, $b$,   is 
 defined below before the inequality (\ref{15b}).
}
\vsn
{\bf Proof}$\;\;$
Let $M_0=0$ and 
 $M_{i+1}= \exp(L h_i)  M_i + \left( 2 D_r/r!  + 1/2 \right)h_i$,
$i=0,1,\ldots, m-1$.  One can check that $M_i\leq M$, $i=0,1,\ldots,m$.  
 We shall show by induction on $i$ that
\be
\sup\limits_{t\in [a,x_i]} \|z(t)-l_{r+1}(t)\| \leq M_i \max\limits_{0\leq j\leq i-1} h_j^r,
\label{prop11}
\ee
where  $\max\limits_{0\leq j\leq -1} =1$. 
\f
For $i=0$, (\ref{prop11}) holds true.  Let (\ref{prop11}) hold for some $i$. Consider the interval $[x_i,x_{i+1}]$.  We have that
\be
\|z(t)-z_i(t)\| \leq \exp(Lh_i)\, \|z(x_i)-y_i\|, \;\;\; t\in [x_i,x_{i+1}].
\label{11}
\ee
We shall now study the local error in $[x_i,x_{i+1}]$ given by 
 $e_{i,j}=\sup\limits_{t\in [x_i,x_{i+1}]} \|z_i(t)-l_{i,j}(t)\|$.  Denoting $H_i(t)=f(t,z_i(t))\;(=z_i'(t))$, we  
let $\bar q_i$ be the Lagrange interpolation polynomial for  $H_i$,
\be
\bar q_i(t) = \sum\limits_{k=0}^{r-1} f(t_k, z_i(t_k))) \prod\limits_{p=0,p\ne k}^{r-1} \frac{ t-t_p}{ t_k-t_p}, \;\;\;\; t\in [x_i,x_{i+1}].
\label{barq}
\ee
Since 
$$
z_i(t)=y_i+\int\limits_{x_i}^t f(\xi,z_i(\xi))\, d\xi,
$$
we have that
\be
\|z_i(t)-l_{i,j+1}(t)\| \leq \int\limits_{x_i}^t \|f(\xi,z_i(\xi)) - \bar q_i(\xi)\|\, d\xi + \int\limits_{x_i}^t \| \bar q_i(\xi)- q_{i,j}(\xi)\|\, d\xi.
\label{12}
\ee
By the Lagrange interpolation error formula applied component by component, we have that
$$
\|f(\xi,z_i(\xi)) - \bar q_i(\xi)\| \leq \frac{1}{r!}\sup\limits_{\alpha\in [x_i,x_{i+1}]}\|H_i^{(r)}(\alpha) \| \prod\limits_{k=0}^{r-1}|\xi-t_k|,
\;\;\; \xi\in [x_i,x_{i+1}].
$$
Furthermore,
$$
\| \bar q_i(\xi)- q_{i,j}(\xi)\| \leq \sum\limits_{k=0}^{r-1} \|f(t_k, z_i(t_k)) - f(t_k, l_{i,j}(t_k)) \|\prod\limits_{p=0,p\ne k}^{r-1} \left|\frac{ \xi-t_p}{ t_k-t_p}\right|,
$$
which yields 
that
$$
\| \bar q_i(\xi)- q_{i,j}(\xi)\| \leq L 
\hat C_r \sup\limits_{t\in [x_i,x_{i+1}]} \|z_i(t)-l_{i,j}(t)\|, \;\;\;\; \xi\in [x_i,x_{i+1}],
$$
where 
$\hat C_r$  only depends on $r$. 
From (\ref{12}) we get for $j=0,1,\ldots$
\be
e_{i,j+1} \leq \frac{1}{r!}\sup\limits_{\alpha\in [x_i,x_{i+1}]}\|H_i^{(r)}(\alpha) \| \,h_i^{r+1} + h_i L 
\hat C_r e_{i,j}.
\label{13}
\ee
By solving (\ref{13}) we get for $j=0,1,\ldots $
$$
e_{i,j}\leq  \frac{1}{r!}\sup\limits_{\alpha\in [x_i,x_{i+1}]}\|H_i^{(r)}(\alpha) \| \,h_i^{r+1} \; \frac{1-(h_iL
\hat C_r) ^j}{1-h_iL
\hat C_r} +
(h_iL
\hat C_r)^j e_{i,0}.
$$
Since 
$$
e_{i,0}\leq h_i  \sup\limits_{\alpha\in [x_i,x_{i+1}]} \|f(\alpha, z_i(\alpha))\|,
$$
for $m$ sufficiently large  (such that $h_iL
\hat C_r \leq 1/2$) we have that  
\be
e_{i,j}\leq  \frac{2}{r!}\sup\limits_{\alpha\in [x_i,x_{i+1}]}\|H_i^{(r)}(\alpha) \| \,h_i^{r+1}  +
  h_i^{j+1}   (L
\hat C_r)^j   \sup\limits_{\alpha\in [x_i,x_{i+1}]} \|H_i(\alpha) \|.
\label{14}
\ee
Note that $H_i^{(r)}(\alpha)$ and $H_i(\alpha)$ can be expressed in terms of partial derivatives of the function $f$
of order $0,1,\ldots, r$, evaluated 
at $(\alpha, z_i(\alpha))$. Due to (\ref{11}) and the inductive assumption, for sufficiently large $m$ we have that
$\|z(\alpha)-z_i(\alpha)\|\leq 1$ for $\alpha\in [x_i,x_{i+1}]$. Hence, 
\be
(\alpha,z_i(\alpha))\in [a,b]\times D, 
\label{15a}
\ee
where the set $D$ is given in (\ref{D}).
This  yields that $\sup\limits_{\alpha\in [x_i,x_{i+1}]}\|H_i^{(r)}(\alpha) \|$ and $\sup\limits_{\alpha\in [x_i,x_{i+1}]}\|H_i(\alpha) \|$
are bounded from above independently of $i$ and $m$ by some numbers $D_r$ and $D_0$, respectively. 
For the final approximation $l_{i,r+1}$ we have from (\ref{14}) that
\be
e_{i,r+1}\leq \frac{2}{r!}  D_r  h_i^{r+1}  + h_i^{r+2}   (L
\hat C_r)^{r+1}  D_0,
\label{15b}
\ee
which yields for sufficiently large $m$ that
\be
e_{i,r+1}\leq   \left( \frac{2}{r!} D_r  + \frac{1}{2} \right)  h_i^{r+1} . 
\label{15c}
\ee
For the final global approximation $l_{r+1}$,  we have for $t\in [x_i,x_{i+1}]$
\be
\|z(t)-l_{r+1}(t)\| \leq \|z(t)-z_i(t)\|  + \|z_i(t)-l_{i, r+1}(t)\| \leq \exp(L h_i) \|z(x_i)-y_i\| + e_{i,r+1}. 
\label{15}
\ee
  By the inductive assumption, we get that
\be
\|z(t)-l_{r+1}(t)\| \leq  \exp(L h_i)  M_i \max\limits_{0\leq j\leq i-1} h_j^r +  \left(\frac{2}{r!}  D_r  +\frac{1}{2} \right)h_i^{r+1}, \;\; t\in [x_i,x_{i+1}].
\label{16}
\ee
Hence,
\be
\sup\limits_{t\in [a,x_{i+1}]}  \|z(t)-l_{r+1}(t)\| \leq M_{i+1} \max\limits_{0\leq j\leq i} h_j^r,
\label{17}
\ee
where $M_{i+1}= \exp(L h_i)  M_i +    \left( 2D_r /r! + 1/2 \right)  h_i$. The induction is finished.
\f
To complete the proof  we recall that $M_i\leq M$ for $i=0,1,\ldots, m,$  where $M$ is given in the statement  of the theorem.\qed
\vsn
Given $m$ and  a mesh $\{x_i\}_{i=0}^m$, the method $\phi^*$  describes the construction of the approximations $l_{i,r+1}$ in $[x_i,x_{i+1}]$
for $i=0,1,\ldots, m-1$. 
Theorem 1 provides the bound on the global error of $\phi^*$.  A selection of the mesh $\{x_i\}_{i=0}^m$ still
remains an open question; we will study this  issue in the next sections.
\vsn
\section { \bf{\Large Local error bounds for $\phi^*$ }}
\vsn
We now extract from the proof of  Theorem 1 bounds on the local error  $e_{i,r+1}$ of the method $\phi^*$. 
By (\ref{14})
\be
e_{i,r+1}\leq \frac{2}{r!}\sup\limits_{\alpha\in [x_i,x_{i+1}]}\|H_i^{(r)}(\alpha) \| \,h_i^{r+1}  +
  h_i^{r+2}   (L
\hat C_r)^{r+1}   \sup\limits_{\alpha\in [x_i,x_{i+1}]} \|H_i(\alpha) \|.
\label{151}
\ee
Let $\beta>0$.  From (\ref{151}),  for sufficiently large $m$ we have 
\be
e_{i,r+1}\leq 2\left(  \frac{1}{r!}\sup\limits_{\alpha\in [x_i,x_{i+1}]}\|H_i^{(r)}(\alpha) \| +\beta \right)\,h_i^{r+1} .
\label{18}
\ee
The function $H_i(\alpha)=f(\alpha,z_i(\alpha))=z_i'(\alpha)$ is not known.
We now show how the term '$\sup$' above can be (asymptotically) replaced by a known quantity.
\f
We take $\bar x_{i+1}$,  $x_i < \bar x_{i+1}\leq b $. For  $\bar h_i=\bar x_{i+1}-x_i$ we assume   that
\be
\bar h_i \leq \gamma   h_i ,
\label{barx}
\ee 
where $\gamma\geq 1 $ is a given number which may depend on $f$, but is independent of $i$ and $m$.
\f
 Let $\bar t_k$, $k=0,1,\ldots, r$
be equidistant points from $[x_i,\bar x_{i+1}]$, $\bar t_0=x_i$, $\bar t_r=\bar x_{i+1}$ (we omit the index $i$).  We construct
an auxillary approximation $\bar l_{i,r+1}$ in the interval $[x_i,\bar x_{i+1}]$ in the same way as we did in (\ref{lagrpol})--(\ref{lost})
in the case of the approximation $l_{i,r+1}$
in the interval $[x_i, x_{i+1}]$, using now as interpolation nodes the points $\bar t_0, \bar t_1,\ldots, \bar t_{r-1}\in [x_i,\bar x_{i+1}]$. 
\f
Let 
$\tilde H_i(\alpha)=f(\alpha,\bar l_{i,r+1}(\alpha))$ and $\tilde H_i[\bar t_0,\bar t_1,\ldots, \bar t_r]$ be the divided difference, computed component by component, for $\tilde H_i$.
  We shall need the bounds stated in the following two lemmas. Recall that $H_i(t)=f(t,z_i(t))$. 
\vsn
{\bf Lemma 1}$\;\;$ {\it Let $f\in F_r$, $\beta>0$ and $\varphi \in (0,1)$. There exists $ m_0$ such that for any $m\geq  m_0$, for any $\{x_i\}_{i=0}^m$ satisfying (\ref{podzial})
and $i=0,1,\ldots, m-1$ we have 
\be
\left( S_i+ \beta\right) (1-\varphi)\leq    
\frac{1}{r!} \sup\limits_{\alpha\in [x_i,x_{i+1}]}\|H_i^{(r)}(\alpha) \|+\beta  \leq  \left( S_i+ \beta\right) (1+\varphi),
\label{lemma1}
\ee
where $S_i = \left\|  H_i[\bar t_0,\bar t_1,\ldots, \bar t_r] \right\|$.
}
\vsn
{\bf Proof}$\;\;$  
Let $\bar \alpha$ be a point from $[x_i, x_{i+1}]$ for which
$$
\sup\limits_{\alpha\in [x_i,x_{i+1}]}\|H_i^{(r)}(\alpha) \|  = \|H_i^{(r)}(\bar \alpha) \| .
$$
For the  $l$th component $H_i^l$ of the function $H_i$, we express the divided difference as
\be
 H_i^l[\bar t_0,\bar t_1,\ldots,\bar t_r] = \frac{1}{r!} \left( H_i^l \right)^{(r)} ( \hat \alpha^l),
\label{pom0}
\ee
where $\hat \alpha^l$ is some  point from $[x_i, \bar x_{i+1}]$.
For any  $l=1,2,\ldots,d$ it holds
$$
\frac{1}{r!} \left| \left(  H_i^l \right)^{(r)} ( \bar \alpha)\right|+\beta =\left(  \frac{1}{r!} \left| \left( H_i^l \right)^{(r)} ( \hat \alpha^l)\right|+\beta\right) (1+\kappa_i^l),
$$
where 
$$
\kappa_i^l= \frac {   \left| \left(H_i^l \right)^{(r)} ( \bar \alpha)\right|/r!  - \left|\left(H_i^l \right)^{(r)} ( \hat \alpha^l)\right|/r! }
{ \left|\left(H_i^l \right)^{(r)} ( \hat \alpha^l)\right|/r! +\beta}.
$$
Similarly to what we have already noticed, the quantity 
$
\left(H_i^l \right)^{(r)} ( t),  \; t\in [x_i, \max\{ x_{i+1}, \bar x_{i+1} \}],
$
can be expressed by values of a continuous function defined by partial derivatives of $f$, evaluated at $(t, z_i(t))$,  where  the 
argument $(t, z_i(t))$ belongs to the compact set $[a,b]\times D$. By the uniform continuity of this function, we have that
$$
\max\limits_{0\leq i\leq m-1} \max\limits_{1\leq l\leq  d} \, \sup\limits_{\bar \alpha, \hat\alpha^l\in [x_i,\max\{ x_{i+1}, \bar x_{i+1} \} ]} 
|\kappa_i^l| \to 0,\;\;\; \mbox { as } m\to \infty.
$$
Hence $|\kappa_i^l| \leq \varphi$ for $m\geq \bar m_0$,   which leads to (\ref{lemma1}).  \qed
\vsn
{\bf Lemma 2}$\;\;$ {\it  Let $f\in F_r$. There exist $\bar C $, $\tilde C$, $ m_0$ such that for any $m\geq  m_0$, for any $\{x_i\}_{i=0}^m$ satisfying (\ref{podzial})
it holds 
\be
\left\| H_i [\bar t_0,\bar t_1,\ldots,\bar t_r]  -   \tilde H_i [\bar t_0,\bar t_1,\ldots,\bar t_r] \right\| \leq \bar C \frac{1}{\bar h_i^r} 
\sup\limits_{t\in [x_i,\bar x_{i+1}]} \| z_i(t)- \bar l_{i,r+1}(t)\| \leq \tilde C \bar h_i ,
\label{lemma2}
\ee 
for $i=0,1,\ldots, m-1$. 
\f
(Here $\bar t_k$ are, as above, the equidistant nodes in $[x_i,\bar x_{i+1}]$; the index $i$ is omitted.)
}
\vsn
{\bf Proof}$\;\;$ The proof follows from the fact that  
$$
H_i [\bar t_0,\bar t_1,\ldots,\bar t_r] - \tilde H_i [\bar t_0,\bar t_1,\ldots,\bar t_r] =\sum\limits_{k=0}^r  
\frac{  f(\bar t_k, z_i(\bar t_k)) -  f(\bar t_k, \bar l_{i,r+1} (\bar t_k))  }
{ \prod\limits_{p=0, p\ne k}^r (\bar t_k-\bar t_p) } ,
$$
and from  (\ref{15c}) applied to the approximation $\bar l_{i,r+1}$ in the interval  $[x_i,\bar x_{i+1}]$ with $h_i$ replaced by $\bar h_i$. \qed  
\vsn
From (\ref{18}),  Lemmas 1 and  2 we get the following computable (asymptotic) upper bound on the local error of the method $\phi^*$.
\vsn
{\bf Theorem 2}$\;\;$ {\it  Let $f\in F_r$,  $\beta>0$ and $ \varphi \in (0,1)$. 
There exists $ m_0$  
such that for all $m\geq  m_0$, for any $\{x_i\}_{i=0}^m$ satisfying (\ref{podzial}), and any $\bar x_{i+1}$ satisfying
(\ref{barx})  it holds 
\be
e_{i,r+1}\leq G_i\,  h_i^{r+1}, \;\;\; i=0,1,\ldots, m-1,
\label{th1}
\ee
where 
\be
G_i=G_i(f)=(8/3) \left( \|\tilde H_i[\bar t_0,\bar t_1,\ldots, \bar t_r] \| +\beta \right)(1+\varphi). 
\label{wspol}
\ee
}
\f
{\bf Proof}$\;\;$ We successively use (\ref{18}), (\ref{lemma1}) and  (\ref{lemma2}). We first get
$$
e_{i,r+1}\leq 2\left( \|  H_i[\bar t_0,\bar t_1,\ldots, \bar t_r] \| +\beta \right)(1+\varphi) \,h_i^{r+1} ,
$$
and next
$$
e_{i,r+1}\leq 2\left( \|\tilde H_i[\bar t_0,\bar t_1,\ldots, \bar t_r] \| + \tilde C \bar h_i +\beta  \right)(1+\varphi) \,h_i^{r+1}  .
$$
For sufficiently large $m$ such that     $\tilde C \bar h_i\leq \beta/3$  we get  (\ref{th1}).    \qed
\vsn
For given $x_i<b$ and $y_i$,  after  selecting a point $\bar x_{i+1}$, we are able to construct
 an auxillary approximation $\bar l_{i,r+1}$ and compute $\tilde H_i[\bar t_0,\bar t_1,\ldots, \bar t_r]$. 
Hence,  for  given $\beta$ and $\varphi$, the coefficient $G_i$  can be effectively computed.  
For further purposes, note that
\be
G_i\geq (8/3)\beta.
\label{pom1}
\ee
On the other hand, due to (\ref{lemma2}), we have in terms of $H_i$ that
$$
G_i\leq (8/3) \left( \| H_i[\bar t_0,\bar t_1,\ldots, \bar t_r] \| +1+\beta \right)(1+\varphi),
$$
for sufficiently large $m$. Due to (\ref{pom0}) and the observation made in the proof of Lemma 1 regarding 
the derivatives of $H_i$, we have the bound
\be
G_i\leq N(f),
\label{pom3}
\ee
where $N(f)$ is independent of $i$ and $m$.
\vsn
\section {\bf{\Large Algorithm with guaranteed local error}} 
\vsn
Let $\e\in (0,1)$.  
Our aim is to select mesh points  $\{x_i^*\}$ in the method $\phi^*$ in order to guarantee that in all steps the local error 
is at most $\e$, 
\be
e_{i,r+1}\leq \e, \;\;\;\;\; i=0,1,\ldots, m-1.
\label{kryt}
\ee
Given $x_i^*$ and $y_i^*$, we set 
\be
\bar x_{i+1}= x_i^*+ \min\left\{ h(\e) , b-x_i^*\right\} ,
\label{aux}
\ee
where $h:(0,1)\to \rr_+$ is any function such that $h(\e)=O(\e^{1/(r+1)})$, as $\e\to 0^+$, with a constant in the $'O'$-notation 
possibly dependent on $f$ (but not on $i$ or the number of subintervals).
 We discuss the form of $h(\e)$ below.
We then compute the auxillary approximation $\bar l_{i,r+1}$ defined in Section 4  before Lemma 1, and  the coefficient $G_i$ from (\ref{wspol}),
with $[x_i,\bar x_{i+1}]$ replaced by $[x_i^*,\bar x_{i+1}]$.
The mesh point $x_{i+1}^*\leq b$ is now selected such that
\be
G_i h_i^{r+1} =\e \;\;\;\; (h_i=x_{i+1}^*-x_i^*), 
\label{Blok}
\ee
that is, we put 
\be
x_{i+1}^*= x_i^*+ \min\left\{ \left( \frac{\e}{G_i} \right)^{1/(r+1)}, b-x_i^*\right\}.
\label{pktsiat}
\ee
With this mesh point $x_{i+1}^*$, the approximation $l_{i,r+1}$ defined in the algorithm $\phi^*$ satisfies, due to Theorem 2, the condition (\ref{kryt}).
\vsn
{\bf Remark 1}$\;\;$ We comment on the choice of $h(\e)$. We can choose 
$\bar x_{i+1}$ 'large' by taking $h(\e)= \e^{1/(r+1)}$, so that
\be
\bar x_{i+1}= x_i^*+ \min\left\{ \e^{1/(r+1)},  b-x_i^*\right\}. 
\label{aux1}
\ee
The condition (\ref{barx}) holds for $\bar x_{i+1}$ and $x_{i+1}^*$, since $x_{i+1}^*-x_i^*= \Theta (\bar x_{i+1}-x_i^*)$, 
 which  follows from (\ref{pom1}) and  (\ref{pom3}).
If we take $h(\e)=O(\e)$, then 
 the condition (\ref{barx}) holds for $\bar x_{i+1}$ and $x_{i+1}^*$, since
$\bar x_{i+1} \leq x_{i+1}^*$  for sufficiently small $\e$.
In this case   $\bar x_{i+1}$  can be arbitrarily close to $x_i^*$.  Our results hold for such $h(\e)$.
 It seems that the faster $h(\e)$ goes to zero with $\e\to 0$, the earlier the asymptotics shows up. 
\vsn
We have arrived at  the following algorithm.
\vsn
{\bf Algorithm ADAPT-MESH}
\vsn
{\bf 1.}$\;\;$ Choose $\e\in (0,1)$, $\beta>0$ and $ \varphi \in (0,1)$. Set $x_0^*=a$ and $y_0^*=\eta$.    \f
{\bf 2.}$\;\;$ Given $x_i^*$ and $y_i^*$, compute $\bar x_{i+1}$ from (\ref{aux}).  Compute the equidistant points
$\bar t_0=x_i^*,\bar t_1,\ldots, \bar t_r=\bar x_{i+1}$  from $[x_i^*,\bar x_{i+1}]$. \f
{\bf 3.}$\;\;$ Compute the approximation $\bar l_{i,r+1}$ in $[x_i^*,\bar x_{i+1}]$ from (\ref{lagrpol}) and (\ref{lposr}), based on
$\bar t_0,\bar t_1,\ldots, \bar t_{r-1} $.  \f
{\bf 4.}$\;\;$ Compute (component by component) the divided difference  $\tilde H_i[\bar t_0,\bar t_1,\ldots, \bar t_r]$, where 
$\tilde H_i(t)=f(t,\bar l_{i,r+1}(t))$. \f
{\bf 5.}$\;\;$   Compute $G_i$ from (\ref{wspol}).  \f
{\bf 6.}$\;\;$  Compute $x_{i+1}^*$ from (\ref{pktsiat}), and the equidistant points from $[x_i^*,x_{i+1}^*]$ with $t_0=x_i^*$ for $r=1$ and 
$ t_0=x_i^*, t_1,\ldots,  t_{r-1}= x_{i+1}^*$ for $r\geq 2$. 
\f
{\bf 7.}$\;\;$  Compute the approximation $ l_{i,r+1}$ in $[x_i^*, x_{i+1}^*]$ from (\ref{lagrpol}) and (\ref{lposr}), 
based on  $ t_0, t_1,\ldots,  t_{r-1}$.  \f
{\bf 8.}$\;\;$  Set $y_{i+1}^*= l_{i,r+1}(x_{i+1}^*)$. If $x_{i+1}^*<b$, then go to {\bf 2} with $i:=i+1$.\f
STOP
\vsn
Steps {\bf 2}--{\bf 5}  can be viewed as a   'prediction' stage, while the final approximation is computed in steps {\bf 6} and {\bf 7}.
\f
We denote by $m^*=m^*(\e)$ the number of subintervals defined by the algorithm ADAPT-MESH.
The local error of the approximation $l_{i,r+1}$ in the interval $[x_i^*,x_{i+1}^*]$ is guaranteed  to be at most $\e$,
$$
\sup\limits_{t\in [x_i^*,x_{i+1}^*]} \| z_i(t) - l_{i,r+1}(t)\| \leq \e, \;\; i=0,1,\ldots, m^*-1,
$$
for sufficiently small $\e$, see (\ref{kryt}).
The cost of each step of the method $\phi^*$, when applied on a given mesh,  
is $\k^*(r)=r^2 +\Theta(r)$  evaluations of $f$ which are needed to produce $l_{i, r+1}$. The algorithm ADAPT-MESH
is additionally equipped with the mesh selection procedure which makes the cost  twice as large.
The  cost of ADAPT-MESH  equals $2r^2+\Theta(r)$  function evaluations per step, which is roughly $2m^*\k^*(r)$ in total. In Section 7 we  compare it
with the cost of other methods  and mesh selection procedures.
\vsn
By the definition of $\{x_i^*\} $,  we have that $x_{i+1}^*-x_i^*=\left(\e/G_i\right)^{1/(r+1)}$, $i=0,1,\ldots, m^*-2$, and 
$x_{m^*-1}^*+\left(\e/G_{m^*-1}\right)^{1/(r+1)}\geq b$. This yields that $m^*$ is the minimal number $m\in \nn$ such that
\be
mS(m)\geq (b-a)\left( \frac{1}{\e}\right) ^{1/(r+1)},
\label{123}
\ee
where 
\be
S(m)=\frac{1}{m} \sum\limits_{i=0}^{m-1}\left(\frac{1}{G_i}\right)^{1/(r+1)} .
\label{124}
\ee
Taking into account the bounds on $G_i$, we have that
\be
\left(\frac{8}{3}\beta \right)^{1/(r+1)} (b-a) \left(\frac{1}{\e}\right)^{1/(r+1)} \leq m^* < N(f)^{1/(r+1)} (b-a) \left(\frac{1}{\e}\right)^{1/(r+1)} +1.
\label{125}
\ee
We see that the mesh selection procedure in the algorithm ADAPT-MESH does not reduce  the speed of growth of the cost 
as $\e\to 0$, with respect to the equidistant mesh.
As in the latter case (see Section 7), we have that
\be
m^*(\e) =\Theta\left( \left(\frac{1}{\e}\right)^{1/(r+1)} \right).
\label{126}
\ee
A potential gain of adaption lies in reducing the coefficient.
\f
Note that the condition (\ref{podzial}) is satisfied for $\{x_i^*\}$. Indeed,   we have using (\ref{125}) that 
\be
\max\limits_{0\leq i\leq m^*-1}(x_{i+1}^*-x_i^*) \leq \left(\frac{3\e}{8\beta}\right)^{1/(r+1)} \leq 
2\cdot 3^{1/(r+1)}\left( \frac{N(f)}{8\beta}\right)^{1/(r+1)}   \frac{b-a}{m^*} ,
\label{127}
\ee
for sufficiently small $\e$. 
Hence, (\ref{podzial}) holds with any  
$$K\geq 2\cdot 3^{1/(r+1)}\left( \frac{N(f)}{8\beta}\right)^{1/(r+1)} \;\; \mbox{ and }\;\; {\ell}(m)\geq  \frac{b-a}{m} .$$
\vsn
\section {\bf{\Large Mesh selection for a general class of methods }} 
\noindent
The mesh selection procedure described above can be applied to a class of methods $\phi$ for solving (\ref{1}),
not only for $\phi^*$.  We assume that for any discretization $\{x_i\}_{i=0}^m$  satisfying (\ref{podzial}), 
 a method $\phi$ successively computes  in each interval $[x_i,x_{i+1}]$  an approximation $l_i$ to $z_i$, 
 with $l_i(x_i)=y_i$ and  $y_{i+1}=l_i(x_{i+1})$,  starting from $x_0=a$, $y_0=\eta$. Global approximation $l$  
 computed by $\phi$ in $[a,b]$ is composed of the approximations $l_i$ in $[x_i,x_{i+1}]$,   $\phi(f)(t)=l(t)=l_i(t)$, $ t\in [x_i,x_{i+1}]$,
$i=0,1,\ldots, m-1$.
We assume that the computation of $l_i$ requires a certain number of evaluations of some functionals on $f$ (information about $f$).
The total number of  evaluations  in a single interval $[x_i,x_{i+1}]$ 
is  $\k_{\phi}(r)$,  where  $\k_{\phi}(r)$  is independent of $i$ and $m$. For instance, for the method  $\phi^*$ the functionals are defined by evaluations of
$f$, and $\k_{\phi^*}(r)= 2r^2+\Theta(r)$.
 We assume about $\phi$ that:
\vsn
{\bf A.}$\;\;$ There are  $\bar\beta, \beta >0$ such that for any $f\in F_r$  there is $m_0$ such that for all $m\geq m_0$, 
for any $\{x_i\}_{i=0}^m$ satisfying (\ref{podzial})
\be
 \sup\limits_{t\in [x_i,x_{i+1}]} \| z_i(t)- l_i(t)\| \leq \bar\beta\left( \frac{1}{r!}\sup\limits_{\alpha\in [x_i,x_{i+1}]} \| z_i^{(r+1)}(\alpha)\| + \beta\right)\, h_i^{r+1},
\;\; i=0,1,\ldots,m-1.
\label{general}
\ee
Assumption {\bf A} 
has been verified for $\phi=\phi^*$ in   (\ref{18}).
Of course, the method $\phi^*$  is not the only example of  $\phi$.  It can also be defined 
in many different ways, e.g.,  by  Taylor's approximation.
\vsn
{\bf Remark 2}$\;\;$  It is easy to see that Theorem 1 (with slightly different constant $M$), 
Lemma~1, Lemma 2 and Theorem 2 hold for $\phi$, with $\beta$ given in assumption {\bf A},  with $l_{i,r+1}$ replaced by $l_{i}$ and 
$\bar l_{i,r+1}$ replaced by $\bar l_{i}$. The coefficient   $G_i$ is now given by
\be
G_i=(4/3) \bar \beta  \left( \|\tilde H_i[\bar t_0,\bar t_1,\ldots, \bar t_r] \| +\beta \right)(1+\varphi).
\label{wspol11}
\ee
We have that
$$
(4/3) \bar\beta \beta \leq G_i \leq N(f),
$$
where $N(f)$ depends on $\bar\beta$, and it  is independent of $i$ and $m$.
\vsn
We now discuss yet  another local error bound for $\phi$. 
We  show that  the local solution $z_i$ in the bound (\ref{general})  can be replaced, at cost of changing a constant, by 
the global solution $z$. 
\vsn
{\bf Lemma 3}$\;\;$ {\it For any $f\in F_r$, $\beta>0$ and $\varphi \in(0,1)$ there is $m_0$ such that for all
$m\geq m_0$, for any $\{x_i\}_{i=0}^m$ satisfying (\ref{podzial}) it holds
\be
\frac{1}{r!} \sup\limits_{\alpha\in [x_i,x_{i+1}]} \| z_i^{(r+1)}(\alpha)\| + \beta = \left( \frac{1}{r!}
\sup\limits_{\alpha\in [x_i,x_{i+1}]} \| z^{(r+1)}(\alpha)\| + \beta\right) \left(1+\kappa_i\right),
\label{lemma3}
\ee
$i=0,1,\ldots, m-1$, for some $\kappa_i$, where $|\kappa_i|\leq \varphi$.
}
\vsn
{\bf Proof}$\;\;$ Let the two sup above be achieved in points $t_1, t_2\in [x_i,x_{i+1}]$, respectively. We have
$$
\frac{1}{r!}\|z_i^{(r+1)}(t_1)\|+ \beta= \left( \frac{1}{r!}\|z^{(r+1)}(t_2)\|+ \beta\right)\left(1+ \kappa_i\right),
$$
where $\kappa_i= (1/r!) \left(\|z_i^{(r+1)}(t_1)\| - \|z^{(r+1)}(t_2)\| \right) /\left( (1/r!)\|z^{(r+1)}(t_2)\| +\beta\right)$.
Hence,
$$
|\kappa_i|\leq \frac{1}{\beta r!}\| z_i^{(r+1)}(t_1)-z^{(r+1)}(t_2)\| \leq \frac{1}{\beta r!}\left( \| z_i^{(r+1)}(t_1)-z^{(r+1)}(t_1)\| +
 \| z^{(r+1)}(t_2)-z^{(r+1)}(t_1)\|\right).
$$
The last term is bounded by $\sup\limits_{|t_1-t_2|\leq h_i} \| z^{(r+1)}(t_1) -z^{(r+1)}(t_2)\|$, so that, due to the uniform continuity of $z^{(r+1)}$, 
it tends to $0$
(uniformly with respect to $i$) as $m\to \infty$. Note that $z_i^{(r+1)}(t_1)$ and $z^{(r+1)}(t_1)$ can be  expressed by partial derivatives of $f$ of order 
$0,1,\ldots, r$ evaluated in $(t_1,z_i(t_1))$ and
$(t_1,z(t_1))$, respectively. Due to Theorem 1 for $\phi$,  we have 
$$
\|z_i(t_1)-z(t_1)\|\leq \exp(Lh_i) \|y_i-z(x_i)\|=O\left(\max\limits_{0\leq j\leq m-1} h_j^r \right).
$$
Hence, $z(t_1),z_i(t_1)\in D$, see (\ref{D}), for sufficiently large $m$.
This and the uniform continuity of the partial derivatives of $f$ in $[a,b]\times D$  yield that the first term also tends to $0$ as $m\to \infty$,
uniformly with respect to $t_1$ and $i$. Hence, $|\kappa_i|$ tends to $0$ as $m\to \infty$,
uniformly with respect to $t_1$, $t_2$  and $i$. This proves the lemma. \qed
\vsn
We now list upper bounds on the local error of $\phi$ that appeared so far. The basic one is given in assumption {\bf A}
\be
\sup\limits_{t\in [x_i,x_{i+1}]} \| z_i(t)- l_i(t)\| \leq c_i \, h_i^{r+1},
\label{general1}
\ee
where 
\be
 c_i= \bar\beta\left( \frac{1}{r!} \sup\limits_{\alpha\in [x_i,x_{i+1}]} \| z_i^{(r+1)}(\alpha)\| + \beta\right) .
\label{coef1}
\ee
The bound (\ref{general1}) involves the local solution $z_i$; it usually appears in the error analysis of a method $\phi$. 
The second one follows from Lemma 3 and  has the form (\ref{general1}) with $c_i$ replaced by
\be
 \bar c_i= \bar\beta\left( \frac{1}{r!} \sup\limits_{\alpha\in [x_i,x_{i+1}]} \| z^{(r+1)}(\alpha)\| + \beta\right)\left(1+\varphi\right) .
\label{coef2}
\ee
We observe that the bounds $\bar c_i\, h_i^{r+1}$, $i=0,1,\ldots, m-1$, only depend on the local behavior of $f$ and
on the mesh $\{x_i\}$. They hold for any method $\phi$ satisfying {\bf A}, and are useful for theoretical reasons. 
Note that 
the function $p(x_i,x_{i+1}) = \bar c_i\, h_i^{r+1}$  is 
an increasing function with respect to $x_{i+1}$ (for fixed $x_i$), and a decreasing function with respect to
$x_i$ (for fixed $x_{i+1}$).  
\f
The third  bound is constructive and it will be used in the algorithm ADAPT-MESH-GEN below.  It is given by (\ref{general1}) with $c_i$ replaced by
$G_i$ from  (\ref{wspol11}).
Note that the coefficients  (\ref{coef2}) and  (\ref{wspol11}) are not overestimated compared to (\ref{coef1}); they are
equivalent up to  a coefficient only dependent on $\varphi$, for sufficiently large $m$. 
It follows from Lemmas 1, 2 and 3 that
for any $m\geq m_0$, any $\{x_i\}_{i=0}^m$ satisfying (\ref{podzial}) and  $i=0,1,\ldots, m-1$   it holds
\be
   \frac{1-\varphi}{1+\varphi}\, \bar c_i  \leq c_i \leq \bar c_i \;\;\; \mbox{ and }\;\;\;   \frac{1-\varphi}{2(1+\varphi)}\, G_i \leq c_i \leq G_i.
\label{twosided1}
\ee
Hence,  for a given $\varphi$,  the three error bounds 
$$
c_i\,h_i^{r+1}, \;\;\;\;  \bar c_i\,h_i^{r+1}, \;\; \mbox{ and } \;\; G_i\,h_i^{r+1}, \;\;\;\;
$$
are equivalent up to  absolute constants (for a fixed $\varphi$). In particular,  they all reflect a local behavior of $f$.
\vsn
As in the case of the method $\phi^*$, we are free to choose the mesh points for $\phi$. 
The following algorithm, very much similar to ADAPT-MESH, describes the mesh selection for $\phi$ that allows us to keep the local error at  level $\e$.
\vsn
{\bf Algorithm ADAPT-MESH-GEN}  
\vsn
{\bf 1.}$\;\;$ Choose $\e\in (0,1)$,  and $ \varphi \in (0,1)$. Set $x_0=a$ and $y_0=\eta$. \f
{\bf 2.}$\;\;$ Given $x_i$ and $y_i$, compute 
\be
\bar x_{i+1}= x_i+ \min\left\{   h(\e)  , b-x_i\right\}. 
\label{aux1.1}
\ee
{\bf 3.}$\;\;$  Compute  an approximation $\bar l_{i}$ in $[x_i,\bar x_{i+1}]$ using $\phi$.\f
{\bf 4.}$\;\;$ For $\tilde H_i(t)=f(t,\bar l_{i}(t))$,
compute (component by component) the divided difference  $\tilde H_i[\bar t_0,\bar t_1,\ldots, \bar t_r]$, where 
  $\bar t_0=x_i,\bar t_1,\ldots, \bar t_r=\bar x_{i+1}$  are the equidistant points from $[x_i,\bar x_{i+1}]$. \f
{\bf 5.}$\;\;$  Compute $G_i= (4/3)\bar\beta \left( \|\tilde H_i[\bar t_0,\bar t_1,\ldots, \bar t_r] \| +\beta \right)(1+\varphi)$.  \f
{\bf 6.}$\;\;$  Compute
\be
x_{i+1}= x_i+ \min\left\{ \left( \frac{\e}{G_i} \right)^{1/(r+1)}, b-x_i\right\}.
\label{pktsiat1}
\ee
\noindent
{\bf 7.}$\;\;$  Compute an approximation $ l_{i}$ in $[x_i, x_{i+1}]$ using $\phi$. \f
{\bf 8.}$\;\;$  Set $y_{i+1}= l_{i}(x_{i+1})$. If $x_{i+1}<b$, then go to {\bf 2} with $i:=i+1$.\f
STOP
\vsn
Since $c_i\,h_i^{r+1}\leq G_i\, h_i^{r+1}\leq \e$, 
 we  see that 
the local error of the approximation $l_{i}$ computed in the step {\bf 7}  is guaranteed  to be at most $\e$,
\be
\sup\limits_{t\in [x_i,x_{i+1}]} \| z_i(t) - l_{i}(t)\| \leq \e, \;\; i=0,1,\ldots, m-1,
\label{33}
\ee
for sufficiently small $\e$.
The cost of each step of ADAPT-MESH-GEN, measured by the number of evaluations of functionals on $f$ needed to produce $l_i$,
is doubled with respect to the cost of $\phi$ applied on a given mesh.
\vsn
\section {\bf{\Large  Adaptive mesh -- the cost analysis}} 
\vsn
 Consider an arbitrary method $\phi$ satisfying {\bf A}, based on a mesh $x_0=a<x_1<\ldots<x_m=b$ for which (\ref{podzial}) holds.
We measure the cost of  $\phi$,  $\cost (\phi,m)$,
by the total number of evaluations of functionals on  $f$ needed for computing $l$ in all subintervals $[x_i,x_{i+1}]$ , i.e.,
\be
\cost (\phi,m) = \k_{\phi}(r)\, m.
\label{koszt}
\ee 
Our goal is to keep the local error at a prescribed level $\e$, see (\ref{33}).  
Assuming that
the goal is achieved by some $\phi$ with some  mesh  $x_0,x_1,\ldots, x_m$, we  
wish to compare $ \cost (\phi,m)$  with the cost of the algorithm ADAPT-MESH.
\vsn
To compare the costs of algorithms,  we use  the local error bound   $c_i\, h_i^{r+1}$ from assumption {\bf A}. 
 We wish to assure that
\be
 c_i\, h_i^{r+1}\leq \e, \;\;\; \; \mbox{ for all } i ,
\label{33.1}
\ee
which implies (\ref{33}).
We define the reference quantity $\hat m(\e)$   as follows.
Let 
\be
\hat k(m)     = \inf \left\{ \max\limits_{0\leq i\leq m-1}\,  c_i\, h_i^{r+1}  : \;\;\; x_0=a\leq x_1\leq \ldots \leq x_{m-1}\leq x_m=b
\;\mbox{ satisfies (\ref{podzial})  } \right\} .
\label{km}
\ee
 Then we define 
\be
\hat m(\e)= \min \left\{ m\in \nn:\;\;  \hat k(m)\leq \e \, \right \}.
\label{me}
\ee
Thus, $\hat m(\e)$ is the minimal number of subintervals $m$ for which there exists  a mesh with $m+1$ points such that
  $ c_i\, h_i^{r+1}\leq \e$, $ i=0,1,\ldots, m-1$.  
\f
 Define similarly to (\ref{km}) and (\ref{me}) the following, technically useful, quantities 
\be
\bar k(m)     = \inf \left\{ \max\limits_{0\leq i\leq m-1}\,  \bar c_i\, h_i^{r+1}  : \;\;\; x_0=a\leq x_1\leq \ldots \leq x_{m-1}\leq x_m=b
\;\mbox{ satisfies (\ref{podzial}) } \right\}
\label{km1}
\ee
and   
\be
\bar m(\e)= \min \left\{ m\in \nn:\;\;  \bar k(m)\leq \e \, \right \}.
\label{me1}
\ee
Since $\bar c_i\, h_i^{r+1}$ is an increasing function of $x_{i+1}$ (for fixed $x_i$) and a decreasing function of $x_i$ (for fixed $x_{i+1}$),
for sufficiently small  $\e>0$ the quantity $\bar m(\e)$ can be computed as follows.
We start with $\hat x_0=a$, and for a given $\hat x_i$, we compute $\hat x_{i+1}$ as the unique solution of  $\bar c_i\, (x_{i+1}-\hat x_i)^{r+1}=\e$. 
Then $\bar m(\e)$ is the minimal $i$ such that $\hat x_i\geq b$.
\vsn
Note further that
for any $m\geq m_0$ and $\{x_i\}_{i=0}^m$ satisfying (\ref{podzial}), it follows from (\ref{twosided1}) that  
$$
 \max\limits_{0\leq i\leq  m-1}\,  c_i\, h_i^{r+1}\leq \max\limits_{0\leq i \leq  m-1}\, \bar c_i\, h_i^{r+1}\leq \frac{1+\varphi}{1-\varphi}\, \max\limits_{0\leq i\leq  m-1}\,  c_i\, h_i^{r+1},
$$
which yields for $m\geq m_0$ that 
\be
 \hat k(m)\leq \bar k(m)\leq  \frac{1+\varphi}{1-\varphi} \hat k(m). 
\label{cc}
\ee
Hence, for any $\e\in (0,\e_0]$  (and fixed $\phi$, $\bar\beta$)
\be
\hat m\left( \e\right) \leq  \bar m(\e)\leq \hat m\left(\frac{1-\varphi}{1+\varphi}\, \e\right).
\label{cc1}
\ee
\vsn
We now compare the cost of ADAPT-MESH with the cost of other algorithms $\phi$ equipped with any mesh selection procedure. 
The number of subintervals computed by ADAPT-MESH is $m^*(\e)$ and the cost
of producing $l_{i,r+1}$ in each subinterval is $\k^*(r)$.  Since the  cost  is doubled  due to the mesh selection,  it holds
\be
\cost(\phi^*,m^*(\e)) = 2\k^*(r) m^*(\e).
\label{kkoszt}
\ee 
The quantities $\hat k(m)$ and $\hat m(\e)$ depend on $\phi$ (and $\bar \beta$). In the following result, we shall use
for clarity the notation $\hat k_{\phi}(m)$ and $\hat m_{\phi}(\e)$. We compare $m^*(\e)$ from ADAPT-MESH (where $\bar \beta=2$) with 
the minimal number of intervals for any other method $\phi$ with $\bar \beta=2$, equipped with any mesh selection strategy.
We have  
\vsn
{\bf Theorem 3}$\;\;$ {\it Let $f\in F_r$, $\varphi\in (0,1)$  and $\phi$ be any method satisfying {\bf A} with $\bar\beta=2$.  Then
there exists $\e_0\in (0,1)$ such that for any $0< \e\leq \e_0$ it holds
\be
\hat m_{\phi} (p_1\e)\leq  \hat m_{\phi^*}(\e) \leq  m^*(\e) \leq  \hat m_{\phi} \left( p \e\right) ,
\label{thm3}
\ee
with $p_1= (1+\varphi)/(1-\varphi)$ and  $p= (1-\varphi)^2/(2(1+\varphi)^2)$.  
\f
Hence,  
\be
\frac{2\k^*(r)}{\k_{\phi}(r)}\,   \cost \left(\phi,  \hat m_{\phi} \left( p_1 \e\right)\right)  \leq  \cost(\phi^*,m^*(\e))  \leq \frac{2\k^*(r)}{\k_{\phi}(r)}\, 
 \cost \left(\phi,  \hat m_{\phi}  \left( p\e \right)\right).
\label{theorem2}
\ee  
}
\vsn
{\bf Proof}$\;\;$ The algorithm ADAPT-MESH   defines  $m^*+1$ points  $x_i^*$ such that for $h_i=x_{i+1}^*-x_i^*$ 
and $G_i$ for $\phi^*$ we have 
\be
G_i\, h_i^{r+1}=\e , \;\; i=0,1,\ldots, m^*(\e)-2,\;\;\;\; \mbox{ and }\;\;\;  G_{m^*-1}\, h_{m^*-1}^{r+1}\leq \e. 
\label{pom10}
\ee
We show the lower bound in (\ref{thm3}). For $\{x_i^*\}_{i=0}^{m^*}$ and $\e$ sufficiently small,  we have that $c_i\,h_i^{r+1}\leq G_i\,h_i^{r+1}\leq \e$
for $i=0,1,\ldots, m^*-1$. This yields that $\hat k_{\phi^*} (m^*)\leq \e$,  which implies that $\hat m_{\phi^*} (\e)\leq m^*(\e)$.
The further  lower bound follows from (\ref{cc1}).
\f
We now show the upper bound. By (\ref{twosided1}), for any $m\geq m_0$, any $\phi$, 
any $\{x_i\}_{i=0}^m$ satisfying (\ref{podzial}) and any $i$, we have for $h_i=x_{i+1}-x_i$ (and $G_i$ for $\phi$) that
$$
G_i\, h_i^{r+1}\leq \frac{2(1+\varphi)}{1-\varphi}\, \bar  c_i\,h_i^{r+1}.
$$
For $\phi=\phi^*$, $m=m^*$ and the mesh $\{x_i^*\}$ given by ADAPT-MESH, due to (\ref{pom10}), 
we have for $h_i=x_{i+1}^*-x_i^*$ and $\e$ sufficiently small that 
$$
\e\leq \frac{2(1+\varphi)}{1-\varphi}\, \bar  c_i\,h_i^{r+1}, \;\;\; i=0,1,\ldots, m^*-2 ,
$$
that is,
\be
\frac{1-\varphi}{2(1+\varphi)}\,\e \leq  \bar  c_i\,h_i^{r+1}, \;\;\; i=0,1,\ldots, m^*-2.
\label{m3}
\ee
The observation made after (\ref{me1})   yields that the points $\hat x_i$ computed for 
the accuracy $(1-\varphi)/(2(1+\varphi))\e$ satisfy in the light of (\ref{m3}) the inequalities 
 $\hat x_i \leq x_i^*$ for $i=0,1,\ldots,m^*-1$. This implies that 
\be
m^*(\e)\leq \bar m\left( \frac{1-\varphi}{2(1+\varphi)}\,\e\right).
\label{uu}
\ee
Finally, (\ref{cc1}) and  (\ref{uu}) give us the desired  inequality (\ref{thm3})
\be
m^*(\e)\leq \hat m_{\phi} \left( \frac{(1-\varphi)^2}{2(1+\varphi)^2}\,\e\right).
\label{uu1}
\ee 
The inequalities for the cost follow immediately.$~ $ \qed
\vsn
Theorem 3 says that the cost of  ADAPT-MESH can only exceed by the constant $2\k^*(r)/\k_{\phi}(r)$ the
cost of any algorithm $\phi$ with $\bar\beta=2$, with any mesh selection strategy such that the local error is
at most $ (1-\varphi)^2/(2(1+\varphi)^2)\,\e$. This accuracy is more demanding  than $\e$. For instance, if we take $\varphi=1/2$, then
$\e$ in the accuracy demand  is replaced by $\e/18$.  
\f 
Observe that adaptive mesh cannot reduce the speed of growth
of the cost as $\e\to 0$. It follows from (\ref{thm3}) and  (\ref{126}) that for the best choice of points we have,
similarly as for the equidistant mesh, that  
$$
\hat m(\e) =\Theta\left( \left(\frac{1}{\e}\right)^{1/(r+1)} \right).
$$
The asymptotics is thus the same. 
Possible  advantage of the adaptive mesh selection is hidden in the size of the quantity $\hat m(\e)$, see the definition
(\ref{me}).  
To  illustrate this, 
we now discuss possible advantage of the algorithm ADAPT-MESH with respect to another algorithm $\phi$ with $\bar\beta=2$, 
based on the uniform mesh. Note that for  any $\phi$,  any $m$ and the uniform mesh we have
$$
  \max\limits_{0\leq i\leq m-1} \bar c_i\, h_i^{r+1} = \bar N(f) \left(\frac{b-a}{m}\right)^{r+1},
$$
where $\bar N(f) =\bar \beta \left( (1/r!)\sup\limits_{t\in [a,b]} \|z^{(r+1)}(t)\| +\beta\right) (1+\varphi)$. Hence,
\be
\bar m^{{\rm equid}}(\e)=\min\left\{m:\;\; \bar N(f) \left(\frac{b-a}{m}\right)^{r+1}\leq \e\right\} = \left\lceil (b-a)\left(\frac{\bar N(f)}
{\e}\right)^{1/(r+1)}  \right\rceil.
\label{nos}
\ee
 It follows from (\ref{twosided1}) that we also have 
$$
  \max\limits_{0\leq i\leq m-1} c_i\, h_i^{r+1} = \Theta\left( \bar N(f) \left(\frac{b-a}{m}\right)^{r+1}\right), \;\;\; \mbox{ as } m\to \infty,
	$$
	and 
$$
\hat m^{{\rm equid}}(\e) = \Theta \left(\left\lceil (b-a)\left(\frac{\bar N(f)}{\e}\right)^{1/(r+1)}  \right\rceil  \right), \;\;\; \mbox{ as } \e\to 0, 
$$
where  constants in the $\Theta$-notation only depend on $\varphi$.
\f
We want to compare $m^*(\e)$ (which decides about the cost of ADAPT-MESH) with  $\hat m^{{\rm equid}}(\e)$
(which decides about the cost of $\phi$ with the equidistant mesh). We have by (\ref{uu}) the following sequence of
inequalities 
$$  
m^*(\e)\leq  \bar m\left( \frac{1-\varphi}{2(1+\varphi)}\,\e\right) \leq \bar m^{{\rm equid}}\left( \frac{1-\varphi}{2(1+\varphi)}\,\e\right)
$$
\be =
   \Theta \left(\left\lceil (b-a)\left(\frac{\bar N(f)}{\e}\right)^{1/(r+1)}  \right\rceil  \right) =\Theta\left( \hat m^{{\rm equid}}(\e)\right),
\label{mpor}
\ee
where again  constants in the $\Theta$-notation only depend on $\varphi$.
\f
The second inequality in (\ref{mpor})  allows us to understand when the adaption pays off.  Just
 below the definition  (\ref{me1}), we gave a comment on how to compute $\bar m(\e)$. 
The comment yields that $\bar m (\e)$ is the minimal number $m\in \nn$ such that
\be
m \bar S(m)\geq (b-a)\left( \frac{1}{\e}\right) ^{1/(r+1)},
\label{123.1}
\ee
where 
\be
\bar S(m)=\frac{1}{m} \sum\limits_{i=0}^{m-1}\left(\frac{1}{\bar c_i}\right)^{1/(r+1)} .
\label{124.1}
\ee
To see this, consult  similar reasoning leading to (\ref{123}) and (\ref{124}). In the case of  $\bar m^{{\rm equid}}\left( \e\right)$ an analogous
condition to (\ref{123.1}) is given in the first equality in (\ref{nos}):
\be
m \frac{1}{\bar N(f)^{1/(r+1)} }\geq (b-a)\left( \frac{1}{\e}\right) ^{1/(r+1)}.
\label{123.2}
\ee
Comparing (\ref{123.1}) and (\ref{123.2}), we see that the second term in (\ref{mpor})  is much less than the third term
if $\bar S(m)^{-1}$ is much smaller  than $\bar N(f)^{1/(r+1)}$. Since 
$$
\bar N(f) =\bar \beta \left( \frac{1}{r!} \sup\limits_{t\in [a,b]} \|z^{(r+1)}(t)\| +\beta\right) (1+\varphi),
$$
and 
$$
\bar c_i= \bar\beta\left( \frac{1}{r!} \sup\limits_{t\in [x_i,x_{i+1}]} \| z^{(r+1)}(t)\| + \beta\right)\left(1+\varphi\right),
$$
we can identify cases  when the gain of adaption is significant. Adaption pays off for functions for which the size of $\|z^{(r+1)}(t)\|$ 
changes significantly  in parts of the interval $[a,b]$. Of course, the second inequality in (\ref{mpor})
 can also turn into equality.
For such functions $f$  there is no gain of adaption.  Translating the above discussion to similar properties of the
cost of the algorithm is straightforward.
\section {\bf{\Large Numerical example }}  
\noindent
We illustrate the performance of the mesh selection mechanism in ADAPT-MESH by an example (other test examples in $ {\rm C}^{++}$ are
in progress), see \cite{wilga}. We consider  a scalar test problem from \cite{NumAlg}  with a parameter $\delta>0$  
\be
z'(t)=\frac{3}{4} (z(t)-1)^{-3/2}, \;\; t\in [0,1],\;\;\; z(0)=1+\delta,
\label{test}
\ee
with the global solution given by $z(t)= \left(\frac{15}{8} t+ \delta^{5/2}\right)^{2/5}+1$.
 The solution with the initial condition $z(x)=y\;\;$ ($x\geq 0$, $y>1$)
is given by 
$$
z_{x,y}(t)= \left(\frac{15}{8} (t-x)+(y-1)^{5/2}\right)^{2/5}+1.
$$
The right-hand side $f(t,y)=\frac{3}{4} (y-1)^{-3/2}$ is a $C^{\infty}$ function for $y>1$.
\f
The problem (\ref{test}) is a typical test problem whose computational difficulty can be controlled by $\delta$; it grows 
 as $\delta$ tends to zero. 
We  use the algorithm ADAPT-MESH with $r=1$, which corresponds to the Euler method equipped with the mesh selection
algorithm, and with $r=2$. 
For the global solution $z$,
we have  that $|z^{(r+1)}(t)|$ for $t$ close to $0$ behaves like $1/\delta^4$ for $r=1$ and  $1/\delta^{6.5}$ for $r=2$.
For $t$ away from $0$,  $|z^{(r+1)}(t)|$  is essentially a constant. That is,  for small $\delta$  we should observe a significant advantage of
adaptive mesh points over the equidistant points.
\f
The computer precision is   $macheps=10^{-15}$. Obviously, since $macheps$ is fixed and computing time is limited, we cannot verify the 
asymptotic behavior of the algorithm as $\e\to 0$; we are only able to see results for some number of values of $\e$. 
\f
Let us now briefly discuss a practical choice of $h(\e)$ in step {\bf 2}  of the algorithm.
  In fixed precision computation, the crucial point is   accuracy of computing 
the divided difference in (\ref{lemma2}) of Lemma 2. 
Due to round off errors in computing both $f(\bar t_k, \bar l_{i,r+1}(\bar t_k))$ and the divided difference, the bound
(\ref{lemma2}) changes to $\tilde C( \bar h_i+ macheps/ \bar h_i^r)$, for some $\tilde C$ dependent on $f$. The minimum 
of the function of $\bar h_i$ is achieved 
for $\bar h_i= (r\cdot macheps)^{1/(r+1)}$. Thus,  
  in step {\bf 2}, neglecting the coefficient dependent on $r$,  we fix $h(\e)$  independently of $\e$ to be $h(\e)= 10^{-15/(r+1)}$.
\f
In step {\bf 5}, we set 
$$
G_i=\left\{ 
\begin{array}{ll}
2 |\tilde H_i[x_i^*, \bar x_{i+1}] |+1 &\;\;\; r=1,\\
4 \left|\tilde H_i[x_i^*, (x_i^*+\bar x_{i+1})/2, \bar x_{i+1}] \right| + 2&\;\;\; r=2.
\end{array}
\right.
$$
(Note that for $r=1$ we have  $l_{i,r+1}=\bar l_{i,r+1}$.) 
The following table shows results computed by ADAPT-MESH for a number of values of  $\delta$ and $\e$. 
We denote 
$$
{\rm MAXERR}=\max\limits_{0\leq i\leq m^*-1} |z_i(x_{i+1}^*)- y_{i+1}^*|,
$$
where $z_i(x_i^*)=y_i^*$, and $m^*$ is the number of intervals computed in ADAPT-MESH.
The value EQUIDIST is the maximal local error of  the  respective  method for $r=1$ or $r=2$ applied on the equidistant 
mesh $x_i=i/m^*$, $i=0,1,\ldots, m^*$,  with the same number of subintervals equal to $m^*$.
In the successive columns we show the values  $m^*$  , MAXERR/$\e$ and  EQUIDIST/$\e$.  
\vsn
{\small
\begin{tabular}{cccccccc}  
$\delta$ & $\e$ & $m^*$  & $\mbox{ MAXERR}/ \e$  & $\mbox{ EQUIDIST}/\e$  & $m^*$  & $\mbox{ MAXERR}/ \e$  & $\mbox{ EQUIDIST}/\e$\\
          & $~$ &  $r=1$  & $r=1$ &  $r=1$ &$r=2$ & $r=2$ &$r=2$\\
 \\$~$\\
0.1 & $10^{-2}$ & 33 & 0.22 & 49.42      & 24&0.03& 26.06\\
     & $10^{-4}$ & 315 &  0.246 & 225.7    &99&0.04&345.62\\
     & $10^{-8}$ &31373 & 0.25 & 424.4       &2081&0.04& 5331.38   \\
     & $10^{-14}$ &31371619 & 0.264 & 428.7    &207780& 0.06& 7409.35    \\
$10^{-2}$  & $10^{-2}$ &41 & 0.22 & $1801.15$       &33&0.04& $1105.64$    \\
     & $10^{-4}$ &390 & 0.25 & $18147.4$              &136& 0.11& $25876.9$\\
     & $10^{-8}$ &38841 & 0.25 & $907049$                 &2821& 0.16& $9.15* 10^6$\\
     & $10^{-14}$ &38839361 & 0.26 & $2.79*10^{6}$         &281583& 0.175& $4.73* 10^{9}$\\
$10^{-3}$  & $10^{-2}$ &43 & 0.22 & $55127.5 $       &32&1.3& $37025.9 $    \\
     & $10^{-4}$ &413 & 0.37 & $5.73*10^5$              &140& 18.65& $8.45*10^5$\\
     & $10^{-8}$ &41109 & 0.49 & $5.6*10^7$                 &2917& 950.194& $4.0* 10^8$\\
   & $10^{-14}$ &41106703 & 0.5 & $1.48*10^{10}$         &291133& 2262.01& $3.35* 10^{12}$\\
$10^{-4}$     & $10^{-2}$ &42 & 1.005 & $1.79*10^{6}$       & 22& 16.14& $1.7*10^{6}$\\
     & $10^{-4}$ &414 & 8.09 & $1.81*10^{7}$                         &121& 336.5& $3.1*10^{7}$\\
     & $10^{-8}$ &41367 & 77.96 & $1.81*10^{9}$                            & 2915&  118505& $1.29*10^{10}$\\
     & $10^{-14}$ &41365164 & 88.96 & $1.71*10^{12}$                        &291276& $7.88*10^7$ & $ 1.28*10^{14}$\\
   
\end{tabular}
} 
\vsn
According to the theory,  for sufficiently small $\e$ the values in the 4th column (for $r=1$) and 7th column (for $r=2$) should be at most 1. 
This is the case for $\delta= 10^{-1}, 10^{-2}$ for $r=1,2$, and $\delta=10^{-3}$ for $r=1$. For small values of $\delta$, the round off errors 
do not allow us to observe the asymptotic behavior of the algorithm, since the value of $\e$ is too large.
Comparison of columns 4 and 5 for $r=1$ and 7 and 8 for $r=2$ shows the gain of the adaptive mesh selection algorithm applied
in ADAPT-MESH  over  the equidistant points. In the test we have computed results for the equidistant mesh with the same number of points.
 We may wish to compare the behavior of adaption with nonadaption using the same number of evaluations of $f$.
For $r=1$, the adaptive method uses 2 function evaluations, while the nonadaptive one only one value. Hence, in this case 
the value in the 5th column should be divided by 4. For $r=2$, the respective numbers are 10 and 4 evaluations, that is
the result in the 8th column should be divided roughly by 16. This does not change the picture --
in both cases, for small $\delta$ the tests show  a  very significant advantage of the adaption over nonadaption. 
\f
We  shortly comment on comparison between the algorithm defined in \cite{NumAlg} 
for scalar autonomous problems and the current algorithm designed for systems of IVPs, for the  test
problem (\ref{test}). 
As it can be  expected,  the algorithm from \cite{NumAlg}    allows us to better treat small
values of $\delta$. This follows from the fact that, roughly speaking,  the step size control in \cite{NumAlg} was based on  two-sided estimates of local errors.
Specific properties of scalar autonomous problems were used in \cite{NumAlg}; they cannot  be extended to systems of initial value problems.
In order  to handle  systems of  IVPs,  the present algorithm  uses  upper local error bounds, see (\ref{th1}) and (\ref{wspol}).    
\vsn
\section {\bf{\Large Conclusions}}  
We have proposed a mesh selection algorithm for systems of IVPs  that (asymptotically)  guarantees a given level of the local
error.  The algorithm only requires evaluations of the right-hand side $f$.
Rigorous analysis of the cost has been given,  including comparison with the best choice of the mesh points,
as well as with the uniform mesh.
\vsn

\end{document}